\documentclass[11pt,twoside]{article}
\usepackage[margin=1.3in]{geometry}

\let\oldbibliography\thebibliography
\renewcommand{\thebibliography}[1]{%
\oldbibliography{#1}%
\setlength{\itemsep}{0pt}%
}

\usepackage[T1]{fontenc}
\usepackage[utf8]{inputenc}
\usepackage{amsmath,amssymb,amsthm,mathtools}
\usepackage{newtxtext,newtxmath,microtype}
\usepackage[hidelinks]{hyperref}
\allowdisplaybreaks[2]
\numberwithin{equation}{section}
\newtheorem{theorem}{Theorem}[section]
\newtheorem{lemma}[theorem]{Lemma}
\newtheorem{proposition}[theorem]{Proposition}

\theoremstyle{remark}
\newcommand{\R}{\mathbb R}

\title{\textbf{Liouville theorems for nonlinear polyharmonic equations with indefinite nonlinearities}}
\author{Juncheng Wei \and Meiqing Xu \and Hui Yang}
\date{\today} 

\begin{document}

\maketitle

\begin{abstract}
We establish a Liouville theorem for nonnegative bounded solutions of the indefinite polyharmonic equation 
\[
(-\Delta)^m u=x_1u^p
\qquad\text{in }\mathbb R^n, 
\]
for all integers $m\ge2$, $n\ge1$, and all $p>1$. The main difficulty lies in extending the moving-plane method to a higher-order equation with a sign-changing coefficient.
We first combine sublinear growth estimates with reflection arguments to establish the super-polyharmonic property in the right half-space, without prescribing boundary conditions on the separating hyperplane.
We then use a doubling and rescaling argument to obtain relative first- and second-derivative estimates on strips, replacing the scalar Harnack estimates used in the second-order case. These estimates allow us to establish monotonicity in the $x_1$-direction within the right half-space $\{x_1>0\}$, and a cutoff argument then yields the desired nonexistence result. 

\medskip 


\medskip 


\end{abstract}

\section{Introduction} 
In this paper, we study nonnegative solutions of the nonlinear polyharmonic equation 
\begin{equation}\label{eq:main}
    (-\Delta)^m u=x_1u^p\qquad\text{in }\R^n,
\end{equation}
where $m\geq1$ and $n\geq1$ are integers, $p>1$, and
$x=(x_1,x')\in\R\times\R^{n-1}$, with $x=x_1$ when $n=1$.
Since the coefficient $x_1$ changes sign across the hyperplane $\{x_1=0\}$, equation \eqref{eq:main} has an indefinite nonlinearity. This linear weight $x_1$ arises naturally in blow-up analysis near regular zeros of a smooth sign-changing coefficient. Indeed, if $K(x_0)=0$ and $\nabla K(x_0)\neq 0$, then, after a suitable rotation of coordinates,
\[
K(x_0+ry)=r\kappa y_1+o(r)
\qquad \text{as } r\to 0,
\]
locally uniformly in $y$, where $\kappa=|\nabla K(x_0)|>0$. Thus, appropriate rescalings of $(-\Delta)^m v=K(x)v^p$ near $x_0$ lead, after normalization, to \eqref{eq:main} as a limiting equation. A Liouville theorem for this model can therefore be used to rule out bounded, nontrivial blow-up profiles near such points and, in turn, to establish a priori estimates for solutions of the original equation.

In the second-order case $m=1$, Berestycki, Capuzzo-Dolcetta and Nirenberg~\cite{BCDN1994} proved that
\begin{equation} 
-\Delta u=x_1u^p \qquad \text{in }\mathbb{R}^n
\label{eq:mainsecd}
\end{equation}
admits no nontrivial nonnegative solutions when $n\geq 3$ and $1<p<(n+2)/(n-1)$. Lin~\cite{Lin1998} established nonexistence of positive entire solutions at the Sobolev critical exponent $p=(n+2)/(n-2)$, while Zhu~\cite{Zhu1999} proved a corresponding result at $p=(n+4)/(n-2)$. In the Sobolev critical case, such Liouville theorems play an important role in deriving a priori estimates for equations prescribing sign-changing scalar curvature; see~\cite{ChenLi1997Scalar,Lin1998}. For the more general class of equations $-\Delta u=h(x_1)f(u)$, Du and Li~\cite{DuLi2005} established Liouville theorems for bounded nonnegative entire solutions. Pol\'a\v{c}ik and Quittner~\cite{PolacikQuittner2005} proved a parabolic Liouville theorem whose stationary case excludes nontrivial bounded nonnegative solutions of \eqref{eq:mainsecd} for every $p>1$. Analogous Liouville theorems have also been established in the fractional setting.
Chen and Zhu~\cite{ChenZhu2016} proved that
$(-\Delta)^{\alpha/2}u=x_1u^p$ in $\mathbb{R}^n$
admits no nontrivial bounded nonnegative solutions for $1\leq\alpha<2$ and $p>1$.
For the full range $0<\alpha<2$, Barrios, Del Pezzo, Garc\'ia-Meli\'an and Quaas~\cite{BarriosEtAl2017} and Chen, Li and Zhu~\cite{ChenLiZhu2019} obtained nonexistence results for bounded positive entire solutions with more general monotone weights and nonlinearities, under their respective structural assumptions.

Unlike the Laplacian and fractional Laplacians of order $0<\alpha<2$, the polyharmonic operator $(-\Delta)^m$ with $m\geq2$ does not, in general, satisfy a maximum principle. 
The super-polyharmonic property thus serves as an important tool in proving Liouville theorems and classification results for higher-order equations; see, for example, \cite{ChenFangLi2013,ChenLi2013,Lin1998Biharmonic,WeiXu1999}. Related sign properties have also been investigated for higher-order fractional equations \cite{AoGonzalezHyderWei2022,CaoDaiQin2021,Xu2023,Yang2026}. These works concern equations with nonnegative right-hand sides. For polyharmonic inequalities with negative right-hand sides, related sign properties were studied by Ng{\^o}~\cite{Ngo2017}.
The present problem \eqref{eq:main} differs from both settings because the nonlinearity $x_1u^p$ changes sign across the hyperplane $\{x_1=0\}$. The Liouville results of \cite{XuYangQCurvature} for the conformal $Q$-curvature equation also cover certain radially symmetric sign-changing coefficients, but require these coefficients to be nonpositive outside a ball. The weight $x_1$ in~\eqref{eq:main} changes sign outside every ball and therefore also falls outside these settings.

The sign change of the coefficient $x_1$ across $\{x_1=0\}$ prevents a direct application of existing super-polyharmonic results for equations with nonnegative right-hand sides to~\eqref{eq:main}.
Although the right-hand side is nonnegative in the half-space $\{x_1>0\}$, the absence of prescribed boundary conditions on $\{x_1=0\}$ presents an additional difficulty. In this paper, we first establish the super-polyharmonic property in this half-space (see Proposition~\ref{prop:base}) and then use it to prove the following Liouville theorem. 

\begin{theorem}\label{thm:main}
Let $m\geq2$ and $n\geq1$ be integers, and let $p>1$. Suppose that
$u\in C^{2m}(\R^n)$ is a nonnegative bounded solution of \eqref{eq:main}. Then $u\equiv0$ in $\R^n$.
\end{theorem}

Our proof is motivated by Lin's moving-plane argument for the second-order equation with the linear weight $x_1$ in \cite{Lin1998}.
Extending this argument to higher order involves two additional difficulties.
First, whereas $-\Delta u\geq0$ in $H=\{x_1>0\}$ follows directly from the second-order equation, the signs of the intermediate Laplacians are not immediate when $m\geq2$.
Moreover, no boundary conditions are prescribed on $\partial H$.
We overcome this difficulty by combining sublinear growth estimates with odd reflection differences and even symmetrization to establish the super-polyharmonic property in $H$.

Second, Lin \cite{Lin1998} uses the scalar Harnack inequality to control normalized translates, a step that does not extend directly to polyharmonic equations.
We instead combine the doubling lemma in \cite{PQS2007} with rescaling to obtain relative first- and second-derivative estimates on each strip $\{0\leq x_1\leq L\}$.
The sign and reflection properties established above are essential in proving that the rescaled limits are constant.
These estimates allow us to continue the moving planes and establish monotonicity in the $x_1$-direction in $H$. Finally, compactly supported cutoffs and repeated integration by parts yield the desired contradiction, replacing Lin's first-eigenvalue comparison and avoiding higher-order boundary terms.

The paper is organized as follows. In Section~\ref{sec:preliminaries}, we establish growth estimates, reflection inequalities, and the super-polyharmonic property in the right half-space, and derive relative derivative estimates on strips. In Section~\ref{sec:proof}, we use the method of moving planes to prove monotonicity in the $x_1$-direction in $H=\{x_1>0\}$ and conclude the proof of Theorem~\ref{thm:main} by a cutoff argument.

\section{Super-polyharmonic Properties and Derivative Estimates}\label{sec:preliminaries}
Set
\[
    H:=\{(x_1,x'):x_1>0\},\qquad
    u_j:=(-\Delta)^j u\quad(0\leq j\leq m),
\]
so that $u_0=u$ and
\begin{equation}
    -\Delta u_j=u_{j+1}\quad(0\leq j<m),
    \qquad u_m=x_1 u^p.
\label{eq:chain}
\end{equation}
Write $M=\|u\|_\infty$.  We may assume $M>0$ and derive a contradiction.
The following maximum principle in an unbounded slab
is classical. We include a direct proof for completeness.

\begin{lemma}[Maximum principle in an unbounded slab]\label{lem:slab}
Let $a<b$ and $S=\{a<x_1<b\}$.  Suppose
$h\in C^2(S)\cap C(\overline S)$ is bounded from below in $S$ and satisfies
\[
    -\Delta h\geq0\quad\text{in }S,
    \qquad h\geq0\quad\text{on }\{x_1=a\}\cup\{x_1=b\}.
\]
Then $h\geq0$ in $S$.
\end{lemma}

\begin{proof}
For $n=1$ this is the ordinary one-dimensional maximum principle.  Assume
$n\geq2$, and set
\[
 \rho(x')=(1+|x'|^2)^{1/2}.
\]
Choose $\alpha>0$ so that $\alpha^2+(n-1)\alpha<\frac{\pi^2}{4(b-a)^2}$, and define
\[
 \Psi(x_1,x')=e^{\alpha\rho(x')}
 \cos\!\left(\frac{\pi}{2(b-a)}\left(x_1-\frac{a+b}{2}\right)\right).
\]
The cosine is positive and bounded below by $1/\sqrt2$ for $a\leq x_1\leq b$. Since
\[
 |\nabla_{x'}\rho|^2=\frac{|x'|^2}{1+|x'|^2}\leq1,
 \qquad
 \Delta_{x'}\rho
 =\frac{n-1}{(1+|x'|^2)^{1/2}}
  -\frac{|x'|^2}{(1+|x'|^2)^{3/2}}\leq n-1,
\]
we have
\[
\begin{aligned}
 -\Delta\Psi &=-\partial_{x_1x_1}\Psi-\Delta_{x'}\Psi =\left[\frac{\pi^2}{4(b-a)^2} -\alpha\Delta_{x'}\rho -\alpha^2|\nabla_{x'}\rho|^2\right]\Psi\\
 &\geq\left[\frac{\pi^2}{4(b-a)^2} -(n-1)\alpha-\alpha^2\right]\Psi>0 \qquad\text{in }S.
\end{aligned}
\]
For every $\varepsilon>0$,
\[
 h(x_1,x')+\varepsilon\Psi(x_1,x')
 \geq\inf_S h+\frac{\varepsilon}{\sqrt2}e^{\alpha\rho(x')}
 \longrightarrow+\infty
 \qquad\text{as }|x'|\to\infty,
\]
uniformly for $a\leq x_1\leq b$, since $\inf_S h> -\infty$.
Moreover, $h+\varepsilon\Psi>0$ on
$\{x_1=a\}\cup\{x_1=b\}$.  If $h+\varepsilon\Psi$ were negative
somewhere, these boundary values and the uniform divergence would force
it to attain a negative minimum at an interior point of $S$.  At that
point, $\Delta(h+\varepsilon\Psi)\geq0$, whereas
\[
 -\Delta(h+\varepsilon\Psi)
 =-\Delta h+\varepsilon(-\Delta\Psi)>0
 \qquad\text{in }S.
\]
This contradiction proves $h+\varepsilon\Psi\geq0$ in $S$.
Letting $\varepsilon\downarrow0$ gives $h\geq0$ in $S$.
\end{proof}

We next estimate the growth of $u_j$ for later use of the maximum principle.

\begin{lemma}\label{lem:growth}
For $0\leq j\leq m-1$ there is a constant $C_j$ such that
\begin{equation}
    |u_j(x_1,x')|\leq C_j(1+|x_1|)^{j/m}
    \qquad\text{in }\R^n.
\label{eq:growth}
\end{equation}
In particular, for every $0\leq j\leq m-1$,
\[
    \lim_{|x_1|\to\infty}
    \sup_{x'\in\mathbb R^{n-1}}
    \frac{|u_j(x_1,x')|}{|x_1|}=0.
\]

\end{lemma}

\begin{proof}
For $0<r\leq1$, apply the unit-ball interior estimate for $(-\Delta)^m$
to the function $z\mapsto u(x+rz)$.  The scaling identities
\[
D_z^{2j}[u(x+rz)]\big|_{z=0}=r^{2j}D^{2j}u(x),\qquad (-\Delta_z)^m[u(x+rz)]=r^{2m}(-\Delta)^mu(x+rz)
\]
give
\begin{equation}
 |D^{2j}u(x)|\leq C\left(
 r^{-2j}\|u\|_{L^\infty(B_r(x))}
 +r^{2m-2j}\|(-\Delta)^mu\|_{L^\infty(B_r(x))}
 \right).
\label{eq:interior-growth}
\end{equation}
For $y\in B_r(x)$ we have $|y_1|\leq|x_1|+r$.
Since $(-\Delta)^m u=x_1u^p$ and
$M=\|u\|_{L^\infty(\mathbb R^n)}$, we have
\[
\|(-\Delta)^m u\|_{L^\infty(B_r(x))}
\le (|x_1|+r)M^p.
\]
 Taking
\[
    r=(1+|x_1|)^{-1/(2m)}
\]
and using $r\leq1$, we obtain
\[
\begin{aligned}
 r^{-2j}\|u\|_{L^\infty(B_r(x))}
 &\leq M(1+|x_1|)^{j/m},\\
 r^{2m-2j}\|(-\Delta)^mu\|_{L^\infty(B_r(x))}
 &\leq M^p(1+|x_1|)^{-1+j/m}(|x_1|+r)\leq M^p(1+|x_1|)^{j/m}.
\end{aligned}
\]
Since $u_j=(-\Delta)^j u$ is a linear combination
of derivatives of order $2j$, it follows from
\eqref{eq:interior-growth} that
\[
|u_j(x)|\le C_j(1+|x_1|)^{j/m},
\]
which proves \eqref{eq:growth}.
The case $j=0$ is the assumed boundedness.
\end{proof}

The following lemma is a standard super-polyharmonic property for bounded functions; cf.~\cite[Corollary~7]{Ngo2017}.
For completeness, we include a proof that also covers the case $n=1$.
 
\begin{lemma}\label{lem:superpoly}
Let $F\in C^{2q}(\R^n)$ be bounded and nonnegative, where $q\geq2$, and
suppose
\[
    (-\Delta)^qF\geq0\qquad\text{in }\R^n.
\]
Then
\begin{equation}
    (-\Delta)^jF\geq0
    \qquad (1\leq j\leq q-1).
\label{eq:superpoly}
\end{equation}
\end{lemma}

\begin{proof}
Write $F_j=(-\Delta)^jF$.  We first prove $F_{q-1}\geq0$.  If
$F_{q-1}(x_0)<0$, then $F_{q-1}$ is superharmonic because
$-\Delta F_{q-1}=F_q\geq0$.  For $0\leq j\leq q$ and $n\geq2$,
define the spherical mean about $x_0$ by
\[
 \overline F_j(r):=\frac{1}{|\partial B_r(x_0)|}
 \int_{\partial B_r(x_0)}F_j(y)\,dS_y\qquad(r>0),
\]
where $dS_y$ is surface measure and $|\partial B_r(x_0)|$ is the surface
area of $\partial B_r(x_0)$.  For $n=1$, the definition is
\[
 \overline F_j(r):=\frac{F_j(x_0-r)+F_j(x_0+r)}{2}
 \qquad(r>0,\ 0\leq j\leq q).
\]
In either case, set $\overline F_j(0)=F_j(x_0)$.
Superharmonicity of $F_{q-1}$ gives
\begin{equation}
    \overline F_{q-1}(r)\leq F_{q-1}(x_0)=-c<0
    \qquad(r\geq0).
\label{eq:negative-mean}
\end{equation}
For every $j<q$ the spherical means satisfy
\begin{equation}
 -\left(r^{n-1}\overline F_j'(r)\right)'
   =r^{n-1}\overline F_{j+1}(r),
 \qquad \overline F_j'(0)=0.
\label{eq:mean-ode}
\end{equation}
Integrating \eqref{eq:mean-ode} gives, for $j<q$,
\[
 \overline F_j(r)
 =F_j(x_0)-\int_0^r s^{1-n}
   \int_0^s \tau^{n-1}\overline F_{j+1}(\tau)\,d\tau\,ds.
\]
In particular, \eqref{eq:negative-mean} implies
\[
 \overline F_{q-2}(r)
 \geq F_{q-2}(x_0)+\frac{c}{2n}r^2.
\]
When $q\geq3$, substituting this bound into the identity for
$\overline F_{q-3}$ gives
\[
 \overline F_{q-3}(r)
 \leq F_{q-3}(x_0)-\frac{F_{q-2}(x_0)}{2n}r^2
      -\frac{c}{8n(n+2)}r^4.
\]
Repeating this integration, we obtain
\begin{equation}
 (-1)^{q}\overline F_0(r)
 \geq c_0 r^{2(q-1)}+O(r^{2q-4})
\label{eq:mean-polynomial}
\end{equation}
for some constant $c_0>0$.
If $q$ is even, \eqref{eq:mean-polynomial} forces
$\overline F_0(r)\to+\infty$, contradicting boundedness.  If $q$ is odd,
it forces $\overline F_0(r)\to-\infty$, contradicting $F\geq0$.
Thus $F_{q-1}\geq0$.

Now $F_{q-2}$ is superharmonic.  If it were negative somewhere, the same
argument, with the chain starting at $F_{q-2}$, would again make the
spherical mean of $F$ either unbounded above or eventually negative.
Hence $F_{q-2}\geq0$.  Descending induction proves
\eqref{eq:superpoly}.
\end{proof}

We next establish the super-polyharmonic property in the right
half-space and prove the reflection inequalities needed below.

\begin{proposition}[Super-polyharmonic property in $H$]\label{prop:base}
For every $0\leq j\leq m-1$, $x_1>0$, and $x'\in\R^{n-1}$,
\begin{equation}\label{eq:base-all}
\begin{aligned}
u_j(x_1,x')&\geq u_j(-x_1,x'),\\
u_j(x_1,x')+u_j(-x_1,x')&\geq0,\\
u_j(x_1,x')&\geq0.
\end{aligned}
\end{equation}
In particular, $u$ satisfies the super-polyharmonic property in $H$.
Moreover, if $u\not\equiv0$, then
\begin{equation}\label{eq:strict-positivity}
u>0\qquad\text{on }\overline H.
\end{equation}

\end{proposition}

\begin{proof}
For $x_1>0$ define the odd reflection differences
\[
    A_j(x_1,x')=u_j(x_1,x')-u_j(-x_1,x')
    \qquad(0\leq j\leq m-1).
\]
They vanish on $x_1=0$ and satisfy
\begin{equation}
 -\Delta A_j=A_{j+1}\quad(0\leq j<m-1),
 \qquad
 -\Delta A_{m-1}
 =x_1\bigl(u(x_1,x')^p+u(-x_1,x')^p\bigr)\geq0.
\label{eq:odd-chain}
\end{equation}
We claim $A_j\geq0$ for all $j$.
We first prove that $A_{m-1}\geq 0$ in $H$.
Fix $T>1$. By Lemma~\ref{lem:growth},
\[
\begin{aligned}
|A_{m-1}(x_1,x')|
\leq |u_{m-1}(x_1,x')|
     +|u_{m-1}(-x_1,x')|\leq 2C_{m-1}(1+T)^{(m-1)/m},
0\leq x_1\leq T.
\end{aligned}
\]
Choose $K\geq 2\max\limits_{0\leq j\leq m-1}C_j$,
independent of $T$, and define
\[
v_T(x_1,x')
:=A_{m-1}(x_1,x')
  +K(1+T)^{(m-1)/m}\frac{x_1}{T}.
\]
Then
\[
v_T(0,x')=0,
\]
and
\[
v_T(T,x')
\geq (K-2C_{m-1})(1+T)^{(m-1)/m}
\geq 0.
\]
Moreover, since the added term is harmonic,
\[
-\Delta v_T=-\Delta A_{m-1}\geq 0
\qquad \text{in } \{0<x_1<T\}.
\]
As $v_T$ is bounded from below,
Lemma~\ref{lem:slab} yields $v_T\geq 0$.
Consequently,
\[
A_{m-1}(x_1,x')
\geq -Kx_1\frac{(1+T)^{(m-1)/m}}{T}.
\]
For each fixed $(x_1,x')\in H$, letting
$T\to\infty$ gives
\[
A_{m-1}(x_1,x')\geq 0,
\]
since $(m-1)/m<1$.

Suppose $A_{j+1}\geq0$.  Then $-\Delta A_j\geq0$.  The same argument,
now adding $K(1+T)^{j/m}x_1/T$, proves $A_j\geq0$.  Descending induction
gives
\begin{equation}
    u_j(x_1,x')\geq u_j(-x_1,x')
    \qquad(0\leq j\leq m-1).
\label{eq:odd-positive}
\end{equation}
Next define the even function
\[
    E_0(x_1,x')=u(x_1,x')+u(-x_1,x').
\]
It is bounded and nonnegative.  For $x_1>0$, \eqref{eq:odd-positive} with
$j=0$ gives
\begin{equation}
 (-\Delta)^mE_0
 =x_1\bigl(u(x_1,x')^p-u(-x_1,x')^p\bigr)\geq0.
\label{eq:even-top}
\end{equation}
The left-hand side  $(-\Delta)^mE_0$ is even in $x_1$, so it is nonnegative in all of
$\R^n$.  Lemma~\ref{lem:superpoly} gives
\begin{equation}
 (-\Delta)^jE_0
 =u_j(x_1,x')+u_j(-x_1,x')\geq0
 \quad(1\leq j\leq m-1).
\label{eq:even-positive}
\end{equation}
Combining \eqref{eq:odd-positive} and \eqref{eq:even-positive},
\[
 u_j(x_1,x')=\frac{ (-\Delta)^jE_0+A_j(x_1,x')}{2}\geq0
 \qquad(x_1>0).
\]

Finally, suppose $u\not\equiv 0$.
Since $-\Delta u=u_1\geq 0$ in $H$, the strong
minimum principle implies that either $u>0$
in $H$ or $u\equiv 0$ there.
The latter, together with the first inequality in \eqref{eq:base-all}, gives $u\equiv0$ globally.
Suppose that $u>0$ in $H$ but
$u(0,x'_0)=0$ for some $x'_0\in\mathbb R^{n-1}$.
Since $u\geq 0$ in $\mathbb R^n$, the function
$t\mapsto u(t,x'_0)$ attains its minimum at $t=0$.
Hence,
\[
\partial_{x_1}u(0,x'_0)=0.
\]
On the other hand, the Hopf boundary lemma
applied in $H$ yields
\[
\partial_{x_1}u(0,x'_0)>0,
\]
a contradiction.
\end{proof}

The next lemma bounds the first and second derivatives of $u$ in terms of $u$ on each strip $0\leq x_1\leq L$.

\begin{lemma}\label{lem:relative}
Assume $u\not\equiv0$.  For every $L>0$ there is $C_L<\infty$ such that
\begin{equation}
    \frac{|\nabla u|}{u}
    +\left(\frac{|D^2u|}{u}\right)^{1/2}
    \leq C_L
    \qquad\text{on }\{0\leq x_1\leq L\}.
\label{eq:relative}
\end{equation}
In particular, after enlarging $C_L$ if necessary,
\begin{equation}
 |\partial_{x_1}\log u|\leq C_L,\qquad
 |\partial_{x_1x_1}u|\leq C_Lu
 \qquad\text{on }\{0\leq x_1\leq L\}.
\label{eq:relative-one-dimensional}
\end{equation}
\end{lemma}

\begin{proof}
Suppose that the assertion fails for some $L>0$.
Define
\[
Q(x)=\frac{|\nabla u(x)|}{u(x)}
+\left(\frac{|D^2u(x)|}{u(x)}\right)^{1/2},
\qquad x\in\overline H.
\]
Since $u>0$ on $\overline H$, the function $Q$
is continuous there.

Set
\[
D=\{x\in\overline H:x_1<L+1\},
\qquad
\Gamma=\{x\in\overline H:x_1\geq L+1\},
\]
and let $\mathcal M(x)=\max\{Q(x),1\}$.
For each $k\geq1$, choose $y_{k,0}$ with
\[
0\leq(y_{k,0})_1\leq L,
\qquad Q(y_{k,0})>4k^2.
\]
Then
\[
\mathcal M(y_{k,0})
\operatorname{dist}(y_{k,0},\Gamma)>2k.
\]
Applying the doubling lemma
\cite[Lemma~5.1]{PQS2007}
with $X=\Sigma=\overline H$, we obtain
$y_k\in D$ such that
\[
\mathcal M(y_k)\geq\mathcal M(y_{k,0})>4k^2, \qquad \operatorname{dist}(y_k,\Gamma)
>\frac{2k}{\mathcal M(y_k)},
\]
and
\[
\mathcal M(x)\leq2\mathcal M(y_k)
\]
whenever $x\in D$ and
$|x-y_k|\leq k/\mathcal M(y_k)$.
Since $\mathcal M(y_k)>1$, we have
$\mathcal M(y_k)=Q(y_k)$. Set
\[
\rho_k=Q(y_k)^{-1}
\leq\frac{1}{4k^2}\longrightarrow0.
\]
The distance estimate ensures that
$B_{k\rho_k}(y_k)\cap\overline H\subset D$.
Consequently,
\begin{equation}\label{eq:doubling}
Q(x)\leq2Q(y_k)
\end{equation}
for every $x\in\overline H$ satisfying
$|x-y_k|\leq k\rho_k$. Define
\begin{equation}
    U_k(z)=\frac{u(y_k+\rho_k z)}{u(y_k)}.
\label{eq:rescale}
\end{equation}
Fix $R>0$ and $k\geq R$. We first consider $z\in B_R$ such that $x:=y_k+\rho_kz=(x_1,x')$ satisfies $x_1\geq0$.
 Since
$|x-y_k|=\rho_k|z|\leq k\rho_k$, the doubling estimate
\eqref{eq:doubling} gives
\[
 |\nabla_z\log U_k(z)|
 =\rho_k\frac{|\nabla u(y_k+\rho_kz)|}{u(y_k+\rho_kz)}
 \leq\rho_k Q(y_k+\rho_kz)
 \leq2\rho_k Q(y_k)=2.
\]
As both $y_k$ and $x$ belong to the convex set $\overline H$,
we have $y_k+s\rho_kz\in\overline H$ for every $s\in[0,1]$.
Since $U_k(0)=1$, integration along this segment yields
\[
 |\log U_k(z)|
 =\left|\int_0^1\nabla_z\log U_k(sz)\cdot z\,ds\right|
 \leq2|z|\leq2R,
\]
and hence
\begin{equation}
 e^{-2R}\leq U_k(z)\leq e^{2R}
 \qquad(z\in B_R,\ y_k+\rho_kz\in\overline H).
\label{eq:rescaled-halfspace-bound}
\end{equation}

Now let $z\in B_R$ correspond to a point
$x=y_k+\rho_kz=(x_1,x')$ with $x_1<0$, and consider its
reflection $(-x_1,x')$ across $\{x_1=0\}$.
Since $y_{k,1}\geq0$, we have
\[
 |(-x_1,x')-y_k| \leq |(x_1,x')-y_k|.
\]
Thus $(-x_1,x')\in B_{\rho_kR}(y_k)\cap H$, so its rescaled coordinate $z^*:=((-x_1,x')-y_k)/\rho_k$ belongs to $B_R$ and satisfies $y_k+\rho_kz^*\in H$. Using \eqref{eq:rescaled-halfspace-bound} at $z^*$ together with \eqref{eq:base-all}, we obtain
\begin{equation}
 U_k(z)=\frac{u(x_1,x')}{u(y_k)}
 \leq\frac{u(-x_1,x')}{u(y_k)}
 =U_k(z^*)\leq e^{2R}.
\label{eq:rescaled-reflection-bound}
\end{equation}
Combining \eqref{eq:rescaled-halfspace-bound} with \eqref{eq:rescaled-reflection-bound}, we conclude that
\[
 0\leq U_k(z)\leq e^{2R}
 \qquad(z\in B_R,\ k\geq R).
\]
The rescaled equation is
\begin{equation}
 (-\Delta)^mU_k(z)
 =\rho_k^{2m}(y_{k,1}+\rho_kz_1)
   u(y_k)^{p-1}U_k(z)^p.
\label{eq:rescaled-equation}
\end{equation}
For every fixed $R>0$, its right-hand side satisfies
\[
 \|(-\Delta)^mU_k\|_{L^\infty(B_R)}
 \leq\rho_k^{2m}(L+1+\rho_kR)M^{p-1}e^{2pR}
 \longrightarrow0.
\]
The interior $W^{2m,2n}$ estimate therefore gives
\[
 \|U_k\|_{W^{2m,2n}(B_R)}
 \leq C_{R,m,n}\left(
 \|U_k\|_{L^\infty(B_{2R})}
 +\|(-\Delta)^mU_k\|_{L^\infty(B_{2R})}\right)
 \leq C_{R,m,n,L,M,p}.
\]
Compact Sobolev embedding on each ball, followed by a diagonal
subsequence over $R=1,2,\ldots$, gives
\begin{equation}
 U_k\longrightarrow U
 \quad\text{in }C^{2m-1,\alpha}_{\mathrm{loc}}(\R^n)
\label{eq:convergence}
\end{equation}
for some $0<\alpha<1/2$.  In addition, for every
$\varphi\in C_c^\infty(\R^n)$,
\[
 \int_{\R^n}U(-\Delta)^m\varphi
 =\lim_{k\to\infty}\int_{\R^n}U_k(-\Delta)^m\varphi
 =\lim_{k\to\infty}\int_{\R^n}(-\Delta)^mU_k\,\varphi=0.
\]
Interior elliptic regularity then makes $U$ smooth, and we have
\begin{equation}
    U\geq0,\qquad U(0)=1,\qquad (-\Delta)^mU=0.
\label{eq:limit}
\end{equation}
The normalization follows from
\[
|\nabla U_k(0)|+|D^2U_k(0)|^{1/2} =\rho_k\frac{|\nabla u(y_k)|}{u(y_k)} +\rho_k\left(\frac{|D^2u(y_k)|}{u(y_k)}\right)^{1/2} =\rho_kQ(y_k)=1.
\]
Passing to the limit in the first and second derivatives gives
\begin{equation}
    |\nabla U(0)|+|D^2U(0)|^{1/2}=1.
\label{eq:normalization}
\end{equation}

It remains to show that $U$ is constant.  In the $z$ coordinates the
plane $x_1=0$ is given by
\[
    z_1=-\frac{y_{k,1}}{\rho_k}\leq0.
\]
Pass to a subsequence such that $-{y_{k,1}}/{\rho_k}\to b\in(-\infty,0]$ or
$-{y_{k,1}}/{\rho_k}\to-\infty$. 
If $-{y_{k,1}}/{\rho_k}\to-\infty$,  then
every compact subset of $\mathbb R^n$
is contained in $\left\{z\in\mathbb R^n:
z_1>-{y_{k,1}}/{\rho_k}\right\}$ for all sufficiently
large $k$.  Proposition~\ref{prop:base} and
the scaling identity
\[
 (-\Delta_z)^jU_k(z)
 =\frac{\rho_k^{2j}}{u(y_k)}u_j(y_k+\rho_kz)\geq0
 \qquad(0\leq j\leq m-1)
\]
then imply, by \eqref{eq:convergence},
\begin{equation}
    (-\Delta)^jU\geq0\quad\text{in }\R^n,\qquad 0\leq j\leq m-1.
\label{eq:limit-chain-positive}
\end{equation}
Now $(-\Delta)^{m-1}U$ is a nonnegative entire harmonic function, hence a
constant $c\geq0$.  If $c>0$, then
$-\Delta\bigl((-\Delta)^{m-2}U\bigr)=c$, and hence
\[
\overline{(-\Delta)^{m-2}U}(r)=\bigl((-\Delta)^{m-2}U\bigr)(0)-\frac{cr^2}{2n},
\]
contradicting $(-\Delta)^{m-2}U\geq0$ for large $r$.  Thus $c=0$.
Descending through the chain shows successively that
$(-\Delta)^jU=0$ for $1\leq j\leq m-1$ and that $U$ is constant.

Suppose instead that $-{y_{k,1}}/{\rho_k}\to b> -\infty$.
Passing to the limit in both the odd and even inequalities in
\eqref{eq:base-all} gives, for $z_1>b$ and $0\leq j\leq m-1$,
\begin{equation}
 \bigl((-\Delta)^jU\bigr)(z)\geq \bigl((-\Delta)^jU\bigr)(2b-z_1,z'),
 \qquad \bigl((-\Delta)^jU\bigr)(z)+\bigl((-\Delta)^jU\bigr)(2b-z_1,z')\geq0.
\label{eq:limit-reflection}
\end{equation}
The function $U(z)+U(2b-z_1,z')$ is entire, even about $z_1=b$, and satisfies
\begin{equation}
    (-\Delta)^j\bigl(U(z)+U(2b-z_1,z')\bigr)\geq0\quad(0\leq j\leq m-1),
    \qquad (-\Delta)^m\bigl(U(z)+U(2b-z_1,z')\bigr)=0.
\label{eq:limit-even-chain}
\end{equation}
The same descending harmonic argument used above proves that
$U(z)+U(2b-z_1,z')\equiv C$ for some constant $C$. Thus $0\leq U\leq C$ on $\R^n$. A bounded entire
$m$-harmonic function is constant: the interior estimate
\[
    |\nabla U(x)|\leq C_mR^{-1}\|U\|_{L^\infty(B_R(x))}
\]
and $R\to\infty$ give $\nabla U(x)=0$.
Thus $U$ is constant in both cases, contradicting
\eqref{eq:normalization}.  This proves \eqref{eq:relative}.
\end{proof}

\section{Proof of Theorem~\ref{thm:main}}\label{sec:proof}

For $\lambda\geq0$, set
\[
 \Sigma_\lambda=\{(x_1,x'):x_1<\lambda\},\qquad
 x^\lambda=(2\lambda-x_1,x'),
\]
and, for $0\leq j\leq m-1$,
\begin{equation}
    w_{j,\lambda}(x)=u_j(x^\lambda)-u_j(x).
\label{eq:w}
\end{equation}
They vanish on $x_1=\lambda$ and satisfy
\begin{equation}
 -\Delta w_{j,\lambda}=w_{j+1,\lambda}
 \quad(0\leq j<m-1),
\label{eq:w-chain}
\end{equation}
while
\begin{equation}
 -\Delta w_{m-1,\lambda}
 =(2\lambda-x_1)u(x^\lambda)^p-x_1u(x)^p.
\label{eq:w-top}
\end{equation}
Denote the right-hand side of \eqref{eq:w-top} by $G_\lambda(x)$ on all of $\Sigma_\lambda$.
When $0<x_1<\lambda$, it can be written as
\begin{equation}
 G_\lambda(x)
 =2(\lambda-x_1)u(x^\lambda)^p+c_\lambda(x)w_{0,\lambda}(x),
\label{eq:G}
\end{equation}
where
\[
c_\lambda(x):=
\begin{cases}
\displaystyle
x_1\frac{u(x^\lambda)^p-u(x)^p}
{u(x^\lambda)-u(x)},
& u(x^\lambda)\neq u(x),\\[2mm]
px_1u(x)^{p-1},
& u(x^\lambda)=u(x).
\end{cases}
\]
Since $p>1$ and $u\geq 0$, we have
$c_\lambda(x)\geq 0$ for $0<x_1<\lambda$.

The following lemma shows that if a reflection difference is negative, then $G_\lambda<0$ at some point with $0<x_1<\lambda$.

\begin{lemma}\label{lem:active}
If, for some $\lambda\geq0$, at least one function
$w_{j,\lambda}$ is negative somewhere in $\Sigma_\lambda$, then there is
$x=(x_1,x')$ with $0<x_1<\lambda$ such that $G_\lambda(x)<0$.  At that point, $0<u(x^\lambda)<u(x)$ and
\begin{equation}
    x_1\left(1-\frac{u(x^\lambda)^p}{u(x)^p}\right)>2(\lambda-x_1)\frac{u(x^\lambda)^p}{u(x)^p}.
\label{eq:active-ineq}
\end{equation}
\end{lemma}

\begin{proof}

Fix $\lambda\geq 0$ and suppose that
$G_\lambda\geq 0$ throughout $\Sigma_\lambda$.
Then
\[
-\Delta w_{m-1,\lambda}=G_\lambda\geq 0.
\] 
For any fixed $\varepsilon>0$ and sufficiently
large $T>0$, define
\[
v_T(x):=w_{m-1,\lambda}(x)
+\varepsilon(\lambda-x_1)
\]
in the slab
\(
S_T:=\{-T<x_1<\lambda\}.
\)
We verify the hypotheses of Lemma~\ref{lem:slab}.

First, by Lemma~\ref{lem:growth}, for
$T>\max\{1,\lambda\}$ and $x\in\overline{S_T}$,
\[
\begin{aligned}
|w_{m-1,\lambda}(x)|
&\leq |u_{m-1}(x^\lambda)|
     +|u_{m-1}(x)|\\
&\leq C_{m-1}\bigl[
(1+2\lambda+T)^{(m-1)/m}
+(1+T)^{(m-1)/m}
\bigr]\\
&\leq C_\lambda(1+T)^{(m-1)/m},
\end{aligned}
\]
where $C_\lambda>0$ is independent of
$T$ and $x'$. 
Since $w_{m-1,\lambda}=0$ on $x_1=\lambda$,
we have $v_T(\lambda,x')=0$.
On the other boundary $x_1=-T$,
the growth estimate gives
\[
v_T(-T,x')
\geq -C_\lambda(1+T)^{(m-1)/m}
     +\varepsilon(\lambda+T)\geq 0
\]
for all sufficiently large $T$, since
\[
\frac{(1+T)^{(m-1)/m}}{\lambda+T}
\longrightarrow 0
\qquad\text{as }T\to\infty.
\]
Moreover, since $\lambda-x_1$ is harmonic,
\[
-\Delta v_T
=-\Delta w_{m-1,\lambda}
=G_\lambda\geq 0
\qquad\text{in }S_T.
\]
Finally, since $\lambda-x_1>0$ in $S_T$,
the growth estimate also yields
\[
v_T(x)\geq w_{m-1,\lambda}(x)
\geq-C_\lambda(1+T)^{(m-1)/m},
\qquad x\in S_T.
\]
Thus $v_T$ is bounded from below in $S_T$,
and all the hypotheses of Lemma~\ref{lem:slab}
are satisfied. Consequently, $v_T\geq 0$
in $S_T$, which implies
\[
w_{m-1,\lambda}(x)
\geq-\varepsilon(\lambda-x_1),
\qquad -T<x_1<\lambda.
\]
For each fixed $x\in\Sigma_\lambda$ and every
$\varepsilon>0$, we may choose $T$ sufficiently
large for the preceding argument and for $-T<x_1$.
Thus the resulting inequality holds at $x$
for every $\varepsilon>0$.
Letting $\varepsilon\downarrow0$ yields
\[
w_{m-1,\lambda}(x)\geq0.
\]

Now suppose that $w_{j+1,\lambda}\geq 0$
in $\Sigma_\lambda$ for some $0\leq j\leq m-2$.
By \eqref{eq:w-chain},
\[
-\Delta w_{j,\lambda}=w_{j+1,\lambda}\geq 0.
\]
Moreover, Lemma~\ref{lem:growth} gives
\[
w_{j,\lambda}(x)\geq
-C_{\lambda,j}(1+T)^{j/m},
\qquad -T<x_1<\lambda.
\]
Since $j/m<1$, the preceding argument,
applied to
$w_{j,\lambda}+\varepsilon(\lambda-x_1)$,
yields $w_{j,\lambda}\geq 0$ in
$\Sigma_\lambda$.
By descending induction,
\[
w_{j,\lambda}\geq 0
\qquad (0\leq j\leq m-1),
\]
contrary to the assumption. 
Therefore $G_\lambda<0$ somewhere.  This cannot occur for $x_1\leq0$, so
$0<x_1<\lambda$.  Both $u(x)$ and $u(x^\lambda)$ are positive by
\eqref{eq:strict-positivity}.  If $u(x^\lambda)\geq u(x)$, then
\[
 G_\lambda=(2\lambda-x_1)u(x^\lambda)^p-x_1u(x)^p\geq2(\lambda-x_1)u(x^\lambda)^p\geq0,
\]
a contradiction.  Thus $u(x^\lambda)<u(x)$.  Dividing
\[
    (2\lambda-x_1)u(x^\lambda)^p-x_1u(x)^p<0
\]
by $u(x)^p>0$ and rearranging gives \eqref{eq:active-ineq}.
\end{proof}

Now we are ready to prove the monotonicity of $u$.

\begin{proposition}[Monotonicity in the $x_1$-direction in $H$]
\label{prop:moving-plane-monotonicity}
Under the hypotheses of Theorem~\ref{thm:main},
\begin{equation}
    u(b,x')\geq u(a,x')
    \qquad(0<a<b,\ x'\in\R^{n-1}).
\label{eq:monotone}
\end{equation}
\end{proposition}

\begin{proof}
Assume that $u\not\equiv0$.

\noindent\textit{Step 1. The base plane and the critical position.} By Proposition~\ref{prop:base}, for $x\in\Sigma_0$ and $0\leq j<m$,
\[
 w_{j,0}(x)=u_j(-x_1,x')-u_j(x_1,x')\geq0.
\]
We can therefore define the critical position by
\begin{equation}
 \lambda_*:=\sup\bigl\{\lambda\geq0:
 w_{j,\mu}\geq0\text{ in }\Sigma_\mu
 \text{ for every }0\leq\mu\leq\lambda, 0\leq j<m\bigr\}.
\label{eq:critical}
\end{equation}

\noindent\textit{Step 2. Proving $\lambda_*=\infty$ and deducing monotonicity.} Suppose $\lambda_*<\infty$. By \eqref{eq:critical}, the comparisons hold
at every $0\leq\mu<\lambda_*$. If $\lambda_*>0$, each fixed
$x\in\Sigma_{\lambda_*}$ belongs to $\Sigma_\mu$ for $\mu$ sufficiently
close to $\lambda_*$ from below. Letting $\mu\uparrow\lambda_*$ and using
continuity gives
\[
 w_{j,\lambda_*}(x)\geq0\qquad(x\in\Sigma_{\lambda_*},\ 0\leq j<m).
\]
If $\lambda_*=0$, the same conclusion follows from Step~1.

For every $\delta>0$, at least one comparison fails at some plane in
$(\lambda_*,\lambda_*+\delta)$. Otherwise $\lambda_*+\delta/2$ would
belong to the set in \eqref{eq:critical}. We may thus choose
$\lambda_k\downarrow\lambda_*$ and points
$x_k=(x_{k,1},x'_k)$ such that
\[
0<x_{k,1}<\lambda_k,
\qquad G_{\lambda_k}(x_k)<0,
\]
as guaranteed by Lemma~\ref{lem:active}.  In particular, $0<u(x_k^{\lambda_k})<u(x_k)$.
After passing to a subsequence, $x_{k,1}\to\tau\in[0,\lambda_*]$.
For large $k$ we have $\lambda_k\leq\lambda_*+1$, so both $x_{k,1}$ and $2\lambda_k-x_{k,1}$ lie in $[0,2\lambda_*+2]$.
By Lemma~\ref{lem:relative}, the estimates \eqref{eq:relative-one-dimensional}
hold throughout this strip with a constant $C>0$ independent of $k$ and $x'$.

{Case 1. $\tau<\lambda_*$.} Since $x_{k,1}\to\tau$, we have
$x_{k,1}<\lambda_*$ for all sufficiently large $k$.
Thus
\[
 w_{0,\lambda_*}(x_k)=u(x_k^{\lambda_*})-u(x_k)\geq0.
\]
The first estimate in \eqref{eq:relative-one-dimensional} therefore yields
\[
\log \frac{u(x_k^{\lambda_k})}{u(x_k)} =\log\frac{u(x_k^{\lambda_*})}{u(x_k)} +\int_{2\lambda_*-x_{k,1}}^{2\lambda_k-x_{k,1}} \partial_{x_1}\log u(s,x_k')\,ds \geq-2C(\lambda_k-\lambda_*).
\]
Thus
\[\frac{u(x_k^{\lambda_k})}{u(x_k)}\to1,\]
while $\lambda_k-x_{k,1}\to\lambda_*-\tau>0$.  This contradicts
\eqref{eq:active-ineq}, since its two sides satisfy
\[
 x_{k,1}\left(1-\left(\frac{u(x_k^{\lambda_k})}{u(x_k)}\right)^p\right)\longrightarrow0,
 \qquad 2(\lambda_k-x_{k,1})\left(\frac{u(x_k^{\lambda_k})}{u(x_k)}\right)^p\longrightarrow2(\lambda_*-\tau)>0.
\]

{Case 2. $\tau=\lambda_*$.} The comparison at $\lambda_*$ gives, for $s>0$,
\[
 \frac{u(\lambda_*+s,x_k')-u(\lambda_*-s,x_k')}{2s}
 \geq0.
\]
Letting $s\downarrow0$ implies
\begin{equation}
    \partial_{x_1}u(\lambda_*,x_k')\geq0.
\label{eq:normal-derivative}
\end{equation}
Let
\[
I_k=[x_{k,1},2\lambda_k-x_{k,1}],
\qquad
J_k=[\min\{\lambda_*,x_{k,1}\},
       2\lambda_k-x_{k,1}].
\]
Set
\[
d_k=\lambda_k-x_{k,1},
\qquad
\delta_k=\lambda_k-\lambda_*.
\]
Since $\tau=\lambda_*$, we have $d_k+\delta_k\to0$.
For every $s\in J_k$, the endpoints of $J_k$ give
$|s-\lambda_*|\leq d_k+\delta_k$. Thus
\[
 |s-x_{k,1}|\leq |s-\lambda_*|+|\lambda_*-x_{k,1}|
 \leq d_k+\delta_k+|d_k-\delta_k|\leq 2(d_k+\delta_k).
\]
By the first estimate in \eqref{eq:relative-one-dimensional},
\[
\begin{aligned}
\left|\log\frac{u(s,x'_k)}{u(x_k)}\right|
&=\left|\int_{x_{k,1}}^s
\partial_{x_1}\log u(\xi,x'_k)\,d\xi\right|\leq C|s-x_{k,1}|\leq 2C(d_k+\delta_k).
\end{aligned}
\]
Consequently, for all sufficiently large $k$,
\[
u(s,x'_k)
\leq e^{2C}u(x_k),
\qquad s\in J_k.
\]
The second estimate in
\eqref{eq:relative-one-dimensional} then gives
\[
|\partial_{x_1x_1}u(s,x'_k)|
\leq C'u(x_k),
\qquad s\in J_k,
\]
where $C'>0$ is independent of $k$. After enlarging $C$, \eqref{eq:normal-derivative} therefore gives
\[
\partial_{x_1}u(s,x_k') =\partial_{x_1}u(\lambda_*,x_k')+\int_{\lambda_*}^s \partial_{x_1x_1}u(\xi,x_k')\,d\xi \geq-Cu(x_k)|s-\lambda_*| \qquad(s\in I_k).
\]
Consequently,
\begin{equation}
\begin{aligned}
 1-\frac{u(x_k^{\lambda_k})}{u(x_k)} &=\frac{u(x_k)-u(x_k^{\lambda_k})}{u(x_k)} =-\frac{1}{u(x_k)} \int_{x_{k,1}}^{2\lambda_k-x_{k,1}}\partial_{x_1}u(s,x_k')\,ds\\
 &\leq C\int_{x_{k,1}}^{2\lambda_k-x_{k,1}}|s-\lambda_*|\,ds \leq 2C (\lambda_k-x_{k,1})(2\lambda_k-\lambda_*-x_{k,1})\longrightarrow0.
\end{aligned}
\label{eq:quadratic-contact}
\end{equation}
Since $0<u(x_k^{\lambda_k})<u(x_k)$, we have
\[
 1-\left(\frac{u(x_k^{\lambda_k})}{u(x_k)}\right)^p=p\int_{\frac{u(x_k^{\lambda_k})}{u(x_k)}}^1s^{p-1}\,ds\leq p\left(1-\frac{u(x_k^{\lambda_k})}{u(x_k)}\right).
\]
Using \eqref{eq:active-ineq} and \eqref{eq:quadratic-contact}, we obtain
\[
\begin{aligned}
 2(\lambda_k-x_{k,1})\left(\frac{u(x_k^{\lambda_k})}{u(x_k)}\right)^p
 &<x_{k,1}\left(1-\left(\frac{u(x_k^{\lambda_k})}{u(x_k)}\right)^p\right)\\
 &\leq px_{k,1}\left(1-\frac{u(x_k^{\lambda_k})}{u(x_k)}\right)
 \leq2pC x_{k,1}(\lambda_k-x_{k,1})(2\lambda_k-\lambda_*-x_{k,1}).
\end{aligned}
\]
Thus
\[
 2\left(\frac{u(x_k^{\lambda_k})}{u(x_k)}\right)^p<2pC x_{k,1}(2\lambda_k-\lambda_*-x_{k,1})\longrightarrow0.
\]
By \eqref{eq:quadratic-contact}, the left-hand side tends to $2$.  This is a contradiction.  Hence $\lambda_*=\infty$.
Since $\lambda_*=\infty$, \eqref{eq:critical} with $j=0$ gives
\begin{equation}
 u(x^\lambda)\geq u(x)
 \qquad(\lambda\geq0, x_1<\lambda).
\label{eq:all-planes}
\end{equation}
Given $0<a<b$, choose $x=(a,x')$ and $\lambda=(a+b)/2$ in
\eqref{eq:all-planes} to obtain \eqref{eq:monotone}.
\end{proof}

\begin{proof}[Proof of Theorem~\ref{thm:main}]
Suppose $u\not\equiv0$. Proposition~\ref{prop:moving-plane-monotonicity}
gives \eqref{eq:monotone}. By \eqref{eq:strict-positivity},
choose $x_{1,0}>0$, $x'_0\in\R^{n-1}$, $r>0$, and $c>0$ such that
\begin{equation}
    u(x_{1,0},x')\geq c\qquad\text{when }|x'-x'_0|\leq r.
\label{eq:cylinder-base}
\end{equation}
For $n=1$, this simply means $u(x_{1,0})\geq c$.  Monotonicity gives
\begin{equation}
    u(x_1,x')\geq c
    \qquad(x_1\geq x_{1,0}, ~ |x'-x'_0|\leq r).
\label{eq:cylinder}
\end{equation} 
Let $X_T=(T,x'_0)$ and take $T\geq\max\{x_{1,0}+r,2r\}$.
For every $x\in B_r(X_T)$,
\[
 |x'-x'_0|\leq|x-X_T|<r,
 \qquad x_1\geq T-r\geq x_{1,0}>0.
\]
Thus $B_r(X_T)\subset H$ and \eqref{eq:cylinder} holds throughout that
ball.  Choose the radii explicitly as
\[
 r_j=r\left(1-\frac{j}{2m}\right)\quad(0\leq j\leq m).
\]
For each $0\leq j<m$ take $\eta_j\in C_c^\infty(B_{r_j}(X_T))$ such that
\[
 0\leq\eta_j\leq1,\qquad
 \eta_j=1\text{ on }B_{r_{j+1}}(X_T),\qquad
 |\Delta\eta_j|\leq C_mr^{-2}.
\]
By Proposition~\ref{prop:base}, $u_j\geq0$ in $H$
for $0\leq j\leq m-1$.
Moreover, $u_m=x_1u^p\geq0$ in $H$.
Thus,
\begin{equation}
\begin{aligned}
 \int_{B_{r_{j+1}}(X_T)}u_{j+1}
 &\leq\int_{B_{r_j}(X_T)}u_{j+1}\eta_j
   =\int_{B_{r_j}(X_T)}(-\Delta u_j)\eta_j
   =\int_{B_{r_j}(X_T)}u_j(-\Delta\eta_j)\\
& \leq C_mr^{-2}\int_{B_{r_j}(X_T)}u_j.
\end{aligned}
\label{eq:one-iteration}
\end{equation}
 Applying \eqref{eq:one-iteration} successively for
$j=m-1,m-2,\ldots,0$ and using $u\leq M$ gives
\begin{equation}
\begin{aligned}
 \int_{B_{r/2}(X_T)}x_1 u^p &=\int_{B_{r_m}(X_T)}u_m \leq C_mr^{-2}\int_{B_{r_{m-1}}(X_T)}u_{m-1}\\
 &\leq(C_mr^{-2})^2\int_{B_{r_{m-2}}(X_T)}u_{m-2} \leq\cdots\leq(C_mr^{-2})^m\int_{B_r(X_T)}u \leq C_m^m M|B_1|r^{n-2m}.
\end{aligned}
\label{eq:upper-final}
\end{equation}
On the other hand, for $x\in B_{r/2}(X_T)$ we have
$x_1\geq T-r/2\geq T/2$ and $u(x)\geq c$.  Therefore
\begin{equation}
 \int_{B_{r/2}(X_T)}x_1 u^p
 \geq \frac{T}{2}c^p|B_{r/2}|
 =\frac{|B_1|}{2^{n+1}}T c^p r^n.
\label{eq:lower-final}
\end{equation}
Combining \eqref{eq:upper-final} and \eqref{eq:lower-final} gives
\[
 T\leq\frac{2^{n+1}C_m^mM}{c^p r^{2m}}
 \qquad\text{for every }T\geq\max\{x_{1,0}+r,2r\}.
\]
Since $T\geq\max\{x_{1,0}+r,2r\}$ is arbitrary and the right-hand side is independent of $T$, choosing $T$ larger than the right-hand side
yields a contradiction.
Therefore, $u\equiv0$ in $\R^n$. 
\end{proof}

\paragraph{Acknowledgements} This work was completed while H. Yang was visiting The Chinese University of Hong Kong. He would like to thank the Department of Mathematics and the Institute of Mathematical Sciences for their hospitality. The research of J. Wei was supported by National Key R\&D Program of China (Grant No. 2022YFA1005602) and Hong Kong General Research Fund (No. 14303125) ``On Fujita equation in the critical or supercritical regime''. The research of M.~Xu is partially supported by the National Key R\&D Program of China (Grant No. 2025YFA1017600) and NSFC (No. 12526202). The research of H.~Yang is supported by NSFC (No. 12671139 and No. 12301140). The authors used AI-assisted tools, principally ChatGPT/Codex, for brainstorming and manuscript editing. AI output was not treated as a mathematical authority. The authors independently worked and verified all theorem statements, proofs, calculations, constants, citations, and final prose, and take full responsibility for the contents of the paper.

\noindent J.~Wei 

\noindent Department of Mathematics \\ 
The Chinese University of Hong Kong, Shatin, N.T., Hong Kong \\[1mm]
Email: \href{mailto:wei@math.cuhk.edu.hk}{\nolinkurl{wei@math.cuhk.edu.hk}} 

\bigskip 

\noindent M.~Xu 

\noindent School of Mathematical Sciences\\
Fudan University, Shanghai 200433, China \\
Email: \href{mailto:meiqing_xu@fudan.edu.cn}{\nolinkurl{meiqing_xu@fudan.edu.cn}}

\bigskip 

\noindent H.~Yang

\noindent School of Mathematical Sciences\\
Shanghai Jiao Tong University, Shanghai 200240, China \\
Email: \href{mailto:hui-yang@sjtu.edu.cn}{\nolinkurl{hui-yang@sjtu.edu.cn}}


\begin{thebibliography}{99}
\addcontentsline{toc}{section}{\refname}

\bibitem{AoGonzalezHyderWei2022}
W.~Ao, M.~del Mar Gonz\'alez, A.~Hyder and J.~Wei,
Removability of singularities and superharmonicity for some
fractional Laplacian equations,
\emph{Indiana Univ. Math. J.}
\textbf{71} (2022), no.~2, 735--766.

\bibitem{BarriosEtAl2017}
B.~Barrios, L.~Del Pezzo, J.~Garc\'ia-Meli\'an and A.~Quaas,
A Liouville theorem for indefinite fractional diffusion equations
and its application to existence of solutions,
\emph{Discrete Contin. Dyn. Syst.}
\textbf{37} (2017), no.~11, 5731--5746.

\bibitem{BCDN1994}
H.~Berestycki, I.~Capuzzo-Dolcetta and L.~Nirenberg,
Superlinear indefinite elliptic problems and nonlinear Liouville theorems,
\emph{Topol. Methods Nonlinear Anal.}
\textbf{4} (1994), no.~1, 59--78.

\bibitem{CaoDaiQin2021}
D.~Cao, W.~Dai and G.~Qin,
Super poly-harmonic properties, Liouville theorems and classification
of nonnegative solutions to equations involving higher-order
fractional Laplacians,
\emph{Trans. Amer. Math. Soc.}
\textbf{374} (2021), no.~7, 4781--4813.

\bibitem{ChenFangLi2013}
W.~Chen, Y.~Fang and C.~Li,
Super poly-harmonic property of solutions for Navier boundary problems
on a half space,
\emph{J. Funct. Anal.}
\textbf{265} (2013), no.~8, 1522--1555.

\bibitem{ChenLi1997Scalar}
W.~Chen and C.~Li,
A priori estimates for prescribing scalar curvature equations,
\emph{Ann. of Math.} (2)
\textbf{145} (1997), no.~3, 547--564.

\bibitem{ChenLi2013}
W.~Chen and C.~Li,
Super polyharmonic property of solutions for PDE systems and its applications,
\emph{Commun. Pure Appl. Anal.}
\textbf{12} (2013), no.~6, 2497--2514.

\bibitem{ChenLiZhu2019}
W.~Chen, C.~Li and J.~Zhu,
Fractional equations with indefinite nonlinearities,
\emph{Discrete Contin. Dyn. Syst.}
\textbf{39} (2019), no.~3, 1257--1268.

\bibitem{ChenZhu2016}
W.~Chen and J.~Zhu,
Indefinite fractional elliptic problem and Liouville theorems,
\emph{J. Differential Equations}
\textbf{260} (2016), no.~5, 4758--4785.


\bibitem{DuLi2005}
Y.~Du and S.~Li,
Nonlinear Liouville theorems and a priori estimates for indefinite
superlinear elliptic equations,
\emph{Adv. Differential Equations}
\textbf{10} (2005), no.~8, 841--860.

\bibitem{Lin1998Biharmonic}
C.-S.~Lin,
A classification of solutions of a conformally invariant fourth order
equation in $\mathbb R^n$,
\emph{Comment. Math. Helv.}
\textbf{73} (1998), no.~2, 206--231.

\bibitem{Lin1998}
C.-S.~Lin,
On Liouville theorem and apriori estimates for the scalar curvature equations,
\emph{Ann. Scuola Norm. Sup. Pisa Cl. Sci.} (4)
\textbf{27} (1998), no.~1, 107--130.

\bibitem{Ngo2017}
Q. A. Ng{\^o},
On the sub poly-harmonic property for solutions to
$(-\Delta)^p u < 0$ in $\mathbb{R}^n$,
\textit{C. R. Acad. Sci. Paris, Ser. I}
\textbf{355} (2017), no.~5, 526--532.

\bibitem{PolacikQuittner2005}
P.~Pol\'a\v{c}ik and P.~Quittner,
Liouville type theorems and complete blow-up for indefinite superlinear
parabolic equations,
in \emph{Nonlinear Elliptic and Parabolic Problems}
(M.~Chipot and J.~Escher, eds.),
Progress in Nonlinear Differential Equations and Their Applications,
vol.~64, Birkh\"auser, Basel, 2005, pp.~391--402.

\bibitem{PQS2007}
P.~Pol\'a\v{c}ik, P.~Quittner and P.~Souplet,
Singularity and decay estimates in superlinear problems via
Liouville-type theorems. I. Elliptic equations and systems,
\emph{Duke Math. J.}
\textbf{139} (2007), no.~3, 555--579.

\bibitem{WeiXu1999}
J.~Wei and X.~Xu,
Classification of solutions of higher order conformally invariant equations,
\emph{Math. Ann.}
\textbf{313} (1999), no.~2, 207--228.

\bibitem{Xu2023}
M.~Xu,
On a super polyharmonic property of a higher-order fractional Laplacian,
\emph{Acta Math. Sci. Ser. B (Engl. Ed.)}
\textbf{43} (2023), no.~6, 2589--2596.

\bibitem{XuYangQCurvature}
M.~Xu and H.~Yang,
Liouville theorems for conformal $Q$-curvature equations,
arXiv:\allowbreak2602.14072, 2026.

\bibitem{Yang2026}
H.~Yang,
Liouville-type theorems and radial symmetry for the Yamabe and
Hardy--H\'enon equations involving higher-order fractional Laplacians,
\emph{Sci. China Math.} (2026), published online.
\href{https://doi.org/10.1007/s11425-025-2642-8}{doi:10.1007/s11425-025-2642-8}.

\bibitem{Zhu1999}
M.~Zhu,
Liouville theorems on some indefinite equations,
\emph{Proc. Roy. Soc. Edinburgh Sect. A}
\textbf{129} (1999), no.~3, 649--661.

\end{thebibliography}
\end{document}